\documentclass[12pt,a4paper]{amsart}
\usepackage[english]{babel}
\usepackage{amssymb}
\usepackage{graphicx}
\usepackage{geometry}
\usepackage{color}
\usepackage{epsfig}
\usepackage{subfigure}
\usepackage{tikz}
\usetikzlibrary{positioning}
\usepackage{caption}
\usepackage{hyperref}  
\hypersetup{colorlinks=true, linkcolor=blue, anchorcolor=blue, 
citecolor=red, filecolor=blue, menucolor=blue,
urlcolor=blue}

\numberwithin{equation}{section}
\newtheorem{theorem}{Theorem}

\newtheorem{definition}[theorem]{Definition}

\newtheorem{remark}[theorem]{Remark}

\def\XXint#1#2#3{{\setbox0=\hbox{$#1{#2#3}{\int}$ }
\vcenter{\hbox{$#2#3$ }}\kern-.58\wd0}}

\makeatletter
\@namedef{subjclassname@2020}{
Mathematics Subject Classification}
\makeatother

\begin{document}

\begin{center}
       \vspace{1cm}

       {\Large\textbf{A local surjection theorem with continuous inverse in Banach spaces}}

       \vspace{0.8cm}

       Ivar Ekeland and \'Eric S\'er\'e
       \vspace{0.8cm}

       \textit{Dedicated to the memory of our friend Antonio Ambrosetti.}
       \end{center}
        \vspace{1.5cm}
        
       \textbf{Abstract.} In this paper we prove a local surjection theorem with continuous right-inverse for maps between Banach spaces, and we apply it to a class of inversion problems with loss of derivatives.\vspace{1cm}


\address{Ivar Ekeland: CEREMADE, Universit\'e Paris-Dauphine, PSL Research University, CNRS, UMR 7534, Place de Lattre de Tassigny, F-75016 Paris, France}
\email{ekeland@ceremade.dauphine.fr}

\address{\'Eric S\'er\'e: CEREMADE, Universit\'e Paris-Dauphine, PSL Research University, CNRS, UMR 7534, Place de Lattre de Tassigny, F-75016 Paris, France}
\email{sere@ceremade.dauphine.fr}\medskip

\bigskip


\vspace{2cm}

\bigskip

\bigskip

\section{Introduction}

In the recent work \cite{ES} we introduced a new algorithm for solving nonlinear functional equations admitting
a right-invertible linearization, but with an inverse losing derivatives. These equations are of the form $F(u)=v$ with $F(0)=0$, $v$ small and given, $u$ small and unknown. The main difference with the
classical Nash-Moser algorithm (see {\it e.g.} \cite{Hormander,TZ}) was that, instead of using a regularized Newton scheme, we constructed a sequence $(u_n)_n$ of solutions to Galerkin approximations of the ``hard" problem and proved the convergence of $(u_n)_n$ to a solution $u$ of the exact equation. Each $u_n$ was obtained thanks to a topological theorem on the surjectivity of maps between Banach spaces, due to one of us in \cite{IE2}. However, this theorem does not provide the continuous dependence of $u_n$ as a function of $v$. As a consequence, nothing was said in \cite{ES} on the existence of a continuous selection of solutions $u(v)$. Theorem \ref{Thm8} of the present work overcomes this limitation thanks to a variant of the topological argument, stated in Theorem \ref{continuous inverse}.
\medskip

In the sequel, ${\mathcal L}(X,Y)$ is the space of bounded linear operators between Banach spaces $X$ and $Y$; the operator norm on this space is denoted by $\Vert \cdot\Vert_{X,Y}$. We first restate the result of \cite{IE2} below for the reader's convenience:\bigskip

\begin{theorem}\label{Ivar}{\rm{\bf \cite{IE2}}}
\label{thm1} Let $X$ and $Y$ be Banach spaces. Denote by $B$ the open ball of radius $R>0$ around the origin in $X$. Let $f:B\rightarrow Y$ be
continuous and G\^{a}teaux-differentiable, with $f\left(  0\right)  =0$.
Assume that the derivative $Df\left(  x\right)  $ has a right-inverse
$L\left(  x\right)  $, uniformly bounded on the ball $B_R$:
\begin{align*}
&\forall (x,k)\in B\times Y\text{, \ }Df\left(  x\right) L\left(  x\right)   \,  k=k\,. \\
&\sup\left\{ \, ||\, L\left(  x\right)  ||_{Y,X} \ :\ \left\Vert
x\right\Vert_X < R\right\}   <m\,.
\end{align*}

Then, for every $y\in Y$ with $\left\Vert y\right\Vert_Y < Rm^{-1}$ there is
some $x\in B$ satisfying:
\[
f\left(  x\right)  =y\;\text{ and }\;\left\Vert x\right\Vert_X \leq m\left\Vert y\right\Vert_Y\,.
\]

\end{theorem}

We recall that in the standard local inversion theorem, one assumes that $f$ is of class $C^1$, with $Df(0)$ invertible and $y$ small. An explicit bound on $\left\Vert y\right\Vert_Y$ is provided by the classical Newton-Kantorovich invertibility condition (see \cite{Ciarlet}) when $f$ is of class $C^2$. The bound $\left\Vert y\right\Vert_Y < Rm^{-1}$ of Theorem \ref{Ivar} is much less restrictive than the
Newton-Kantorovich condition, at the price of
losing uniqueness, even in the case when $L(x)$ is also a left inverse of $Df(x)$. To illustrate this, we consider a finite-dimensional example.\medskip

{\bf Example A.} We take $X=Y=\mathbb C$ viewed as a $2$-dimensional real vector space and $f(z)=(2+z)^n-2^n\ ,$ for any complex number $z$ in the open disc of center $0$ and
radius $R=1$ (here $n$ is a positive integer). In that case $Df(z)$ is the multiplication by $n(2+z)^{n-1}$ and $L(z)$ is the multiplication by $n^{-1}(2+z)^{1-n}$,
so $f$ satisfies the assumptions of Theorem 1 for $R=1$ and any real number $m>n^{-1}$. Thus, Theorem \ref{thm1} tells us that the equation $f(z)=Z$
has a solution of modulus less than or equal to $m \vert Z \vert$, provided $Z$ has modulus less than $m^{-1}$. However uniqueness does not hold: the solutions of the algebraic equation $(2+z)^n-2^n=Z$ are of the form $z_k=4ie^{i\frac{k\pi}{n}}\sin{\frac{k\pi}{n}}+O(2^{-n})$, and for
 $\vert Z\vert>4\pi$ the three solutions $z_0$, $z_1$, $z_{-1}$ lie in the closed disc $\{\vert z\vert\leq n^{-1}\vert Z\vert\}$ when $n$ is large enough. 
Yet, there is a unique continuous function $g$ such that $g(0)=0$ and $f\circ g(Z)=Z$ for all complex numbers $Z$ of modulus less than $1/m$. This continuous selection is $g(Z)=z_0=2(\,(1+2^{-n}Z)^{\frac{1}{n}}-1)$ with $(\rho e^{it})^{\frac{1}{n}}=\rho^{1/n}e^{it/n}$, $\forall(\rho,t)\in (0,\infty)\times(-\pi,\pi)$.\medskip

This example raises the following question: in the general case, can we select a solution $x$ depending continuously on $y$, even in infinite dimension and when $Df(x)$ does not have a left inverse? The following theorem gives a positive response, under mild additional assumptions:

\begin{theorem}\label{continuous inverse}
Let $X$, $Y$ be two Banach spaces. Denote by $B$ the open ball of radius $R>0$ around the origin in $X$.
Consider a map $f:B\rightarrow Y$ with $f\left(
0\right) =0$. We assume the following:\medskip

\noindent
$(i)$ $f$ is Lipschitz continuous and G\^ateaux-differentiable on $B$.\medskip

\noindent
$(ii)$ There are a function $L:\, B \to 
\mathcal{L}\left( Y,X\right),$ a constant $a<1$ and, for any $(x,w)\in B\times Y$, a positive radius $\alpha(x,w)$ such that,
if $\Vert
x'-x\Vert_X < \alpha(x,w)$  then $x'\in B$ and
\begin{equation*}
\left\Vert \left( Df\left( x'\right) \circ L\left( x\right)-I_Y\right) w \right\Vert_Y \leq a \Vert w \Vert_Y\;. \label{8}
\end{equation*}

\noindent
$(iii)$ There is some $m<\infty $ such that:%
\[
\sup \left\{ \left\Vert L\left( x\right) \right\Vert_{Y,X} \ :\ x\in B\right\} <m\;.
\]

Denote by $B^{\prime }\subset Y$ the open ball of radius $R^{\prime }:=\left(
1-a\right) R m^{-1}$ and center $0$. Then there is a continuous map $g:B^{\prime }\rightarrow
B$ such that:
\begin{equation*}
\forall y\in B^{\prime }\ , \ \quad \Vert g(y)\Vert_X\leq \frac{m}{1-a}\Vert y \Vert_Y\quad\hbox{ and }\quad   f\circ g\left( y\right) =y\;.
\end{equation*}
If, in addition, one has:\medskip

\noindent
$(iv)$ $f$ is Fr\'echet differentiable on $B$, $Df(x)$ has a left-inverse for all $x\in B\,$ and there is a non-decreasing function $\varepsilon:\,(0,\infty)\to (0,\infty)$ with $\lim\limits_{t\to 0}\varepsilon(t)=0$, such that for all $x_1,\,x_2$ in  $B\,$,
\begin{equation*}
\Vert f(x_2)-f(x_1)-Df(x_1)(x_2-x_1)\Vert_Y\leq \varepsilon(\Vert x_2-x_1\Vert_X) \Vert x_2-x_1\Vert_X\,;
\end{equation*}
then $g$ is the unique continuous right-inverse of $f$ defined on $B'$ and mapping $0_Y$ to $0_X$.
\end{theorem}
\begin{remark}\label{local} If a function $f$ satisfies the assumptions $(i)$, $(ii)$ and $(iii)$ then, for every $x_0\in B$, taking the radius $R_{x_0}=R-\Vert x_0\Vert_X$, one can apply Theorem 1 to the function $z\in B_X(0,R_{x_0})\mapsto f(x_0+z)-f(x_0)$, and one concludes that the restriction of $f$ to $B_X(x_0,R_{x_0})$ has a continuous right-inverse $g_{x_0}$ defined on $B_Y(f(x_0),(1-a)R_{x_0}m^{-1})$ and such that $\Vert g_{x_0}(y)-x_0\Vert_X\leq \frac{m}{1-a}\Vert y-f(x_0) \Vert_Y$ for all $y\in B_Y(f(x_0),(1-a)R_{x_0}m^{-1})$. If, in addition, $f$ satisfies $(iv)$, then $g_{x_0}$ is the unique continuous right-inverse of $f$ defined on the ball $B_Y(f(x_0),(1-a)R_{x_0}m^{-1})$ and mapping $f(x_0)$ to $x_0$.
\end{remark}

\begin{remark}\label{inverse} Assumption (ii) implies that $Df(x)$ has a right inverse $\hat{L}$ such that $\Vert \hat{L}\Vert_{Y,X}\leq (1-a)^{-1}\Vert L\Vert_{Y,X}$. Indeed, taking $P= I_Y-Df\left( x\right) \circ L\left( x\right)$ one can
choose $\hat{L}:=L\circ\left(\sum_{k=0}^\infty P^k\right) $.\medskip

Conversely, Assumption (ii) is satisfied, for instance, if (i) and (ii') hold true, with:\medskip

\noindent (ii') For each $x\in B$, $Df(x)$ has a right-inverse $L(x)\in {\mathcal L}(Y,X)\,$. Moreover, the map
$x\to Df(x)$ is continuous for the strong topology of $B$ and the strong operator topology of ${\mathcal L}(X,Y)\,$: in other words, if $\Vert x_n-x\Vert_X \to 0$ then, for any $v\in X$, $\left\Vert \left( Df(x_n)- Df(x)\right) v\right\Vert_Y\to 0$.\end{remark}

The function $f$ of Example A satisfies the assumptions $(i)$, $(ii')$ and $(iii)$. In that finite-dimensional case, $f$ is of course differentiable in the classical sense of Fr\'echet. Let us give an example for which Fr\'echet differentiability does not hold.
\medskip

{\bf Example B.} Let $\phi\in C^1({\bf R},{\bf R})$
with $\phi'$ bounded on ${\bf R}$ and $\inf_{\bf R} \phi' >0\,.$ The Nemitskii operator $\Phi:\,u\in L^p({\bf R})\to \phi\circ u \in L^p({\bf R}) \,,\; 1\leq p<\infty\,,$ is not Fr\'echet differentiable when $\phi'$ is not constant \cite{Skripnik 1,Skripnik 2}. However, $\Phi$ satisfies conditions $(i)$, $(ii')$ and $(iii)$ for any $r>0$ and $m>(\inf_{\bf R} \varphi')^{-1}\,.$
Therefore, Theorem \ref{continuous inverse} applies to $\Phi$, but the inverse $\Psi$ is easily found without the help of this theorem, as a Nemitskii operator: $\Psi(u)=\psi\circ u$ with $\psi=\phi^{-1}$.
\medskip

It turns out that any function $f$ satisfying $(i)$ has the {\it Hadamard differentiability property} which is stronger than the G\^ateaux differentiability and that we recall below:

\begin{definition}
Let $X$ and $Y$ be normed spaces. A map $f:X\rightarrow Y$ is called
Hadamard differentiable at $x$, with derivative $Df\left( x\right) \in 
\mathcal{L}\left( X,Y\right) $, if, for every sequence $v_{n}\rightarrow v$
in $V$ and every sequence $h_{n}\rightarrow 0$ in $\mathbb{R}$, we have:%
\begin{equation}
\lim_{n}\frac{1}{h_{n}}\left( f\left( x+h_{n}v_{n}\right) -f\left( x\right)
\right) =Df\left( x\right) v \;. \label{a}
\end{equation}
\end{definition}

This notion is weaker than Fr\'{e}chet differentiability but in finite dimension Hadamard and Fr\'{e}chet differentiability are equivalent. On the other hand, Hadamard differentiability is stronger than G\^{a}teaux differentiability, but if a map $f$ is G\^{a}teaux differentiable and Lipschitz, then it is Hadamard differentiable (see, e.g., \cite{Fern}). In particular, the functions $f$ of Theorems \ref{Ivar}, \ref{continuous inverse} are Hadamard differentiable.

 Note that the chain rule holds true for Hadamard differentiable functions, while this is not the case with G\^{a}teaux differentiability (see \cite{Fern}). Hadamard differentiable functions are encountered for instance in statistics \cite{Mises,Vaart,Fern} and in the bifurcation theory of nonlinear elliptic partial differential equations \cite{Ev-St}.\medskip

The paper is organized as follows. In Section \ref{proof} we prove Theorem \ref{continuous inverse}. In Section \ref{hard} we state the hard surjection theorem with continuous right-inverse (Theorem \ref{Thm8}) that can be proved using Theorem \ref{continuous inverse} and proceeding as in \cite{ES}. Finally, under additional assumptions we state and prove the uniqueness of the continuous right-inverse (Theorem \ref{Thm9}).

\section{Proof of Theorem \ref{continuous inverse}}\label{proof}

In \cite{IE2}, Theorem \ref{Ivar} was proved by applying Ekeland's variational principle in the Banach space $X$, to the map $x \mapsto \Vert f(x)-y\Vert_Y$. This principle provided the existence of an approximate minimiser $\underbar{x}$. Assuming that $\Vert f(\underbar{x})-y\Vert_Y>0$ and considering the direction of descent $L(\underbar{x})\left(y-f(\underbar{x})\right)$, a contradiction was found.
Thus, $f(\underbar{x})-y$ was necessarily equal to zero and $\underbar{x}$ was the desired solution of the equation $f(x)=y$. However, there was no continuous dependence of $\underbar{x}$ as a function of $y$. In order to obtain such a continuous dependence, it is more convenient to solve {\it all} the equations $f(x)=y$ for all possible values of $y\in B'$ simultaneously, by applying the variational principle in a functional space of continuous maps from $B'$ to $X$.
The drawback is that it is more difficult to construct a direction of descent, as this direction should be a continuous function of $y$. In order to do so, we use an argument inspired of the classical pseudo-gradient construction for $C^1$ functionals in Banach spaces \cite{Palais}, which makes use of the paracompactness property of metric spaces.\medskip

Consider the space $\mathcal{C}$ of continuous maps $g:B^{\prime
}\rightarrow X$ such that $\Vert y \Vert^{-1}g(y)$ is bounded on $\dot{B}'$, with the notation $\,\dot{B}':=B'\setminus\{0\}$. Endowed with the norm
$$\Vert g \Vert_{\mathcal{C}}= \sup\limits_{\dot{B}'}\Vert y \Vert^{-1}\Vert g(y)\Vert\,,$$
$\mathcal C$ is a Banach space. Consider the function:%
\begin{eqnarray*}
\varphi \left( g\right) &:=&\sup_{y\in \dot{B}'}\Vert y \Vert^{-1}\left\Vert f\circ g\left(
y\right)  -y\right\Vert \ \hbox{ if }\ \Vert g \Vert_{\mathcal{C}} \leq\frac{m}{1-a} \\
\varphi \left( g\right) &:=&+\infty \ \hbox{ otherwise. }
\end{eqnarray*}
The function $\varphi$ is lower semi-continuous on $\mathcal{C}$ and its restriction to the closed ball $\{g\in \mathcal{C}\,:\;\Vert g \Vert_{\mathcal{C}}\leq\frac{m}{1-a}\}$ is finite-valued.
In addition, we have:
\begin{eqnarray*}
\varphi \left( 0\right) &=&\sup\limits_{\dot{B}'}\Vert y \Vert^{-1}\Vert f(0)-y\Vert =1 \\
\varphi \left( g\right)&\geq& 0\ ,\ \ \forall \,g \in \mathcal{C }\,.
\end{eqnarray*}

Choose some $m_{0}$ with:
\begin{equation*}
\sup \left\{ \left\Vert L\left( x\right) \right\Vert_{Y,X} \ :\ x\in B\right\}
<m_{0}<m \;.
\end{equation*}

By Ekeland's variational principle \cite{IE1}, there exists some $g_{0}\in \mathcal{C }$ such that:%
\begin{eqnarray}
\varphi \left( g_{0}\right) &\leq &1  \label{1} \\
\left\Vert g_{0}-0\right\Vert_{\mathcal{C}} &\leq &\frac{m_0 }{1-a}  \label{2} \\
\forall g\in \mathcal{C}\,,\  &&\varphi \left( g\right) \geq \varphi \left(
g_{0}\right) -\frac{(1-a)\varphi \left( 0\right) }{m_{0}}\left\Vert
g-g_{0}\right\Vert_{\mathcal{C}} \;. \label{3}
\end{eqnarray}

Equation \eqref{2} implies that $g_{0}$ maps $B^{\prime }$ into the open ball of center $0_X$ and radius $m_{0}(1-a)^{-1}R'=Rm_0m^{-1}<R$, and the last
equation can be rewritten:%
\begin{equation}
\forall g\in \mathcal{C}\,,\  \varphi \left( g\right)  \geq \varphi \left( g_{0}\right) -\frac{1-a}{m_{0}} \left\Vert g-g_0\right\Vert_{\mathcal{C}} \;.
 \label{opti}
\end{equation}

If $\varphi \left( g_{0}\right) =0$, then $f\left( g_{0}\left( y\right)
\right) -y=0$ for all $y\in B^{\prime }$ and the existence proof is over. If not, then 
$\varphi \left( g_{0}\right) >0$ and we shall derive a contradiction.
In order to do so, we are going to build a deformation $g_t$
of $g_0$ which contradicts the optimality property (\ref{3}) of $g_0$.\medskip

Let $a<a'<1$ be such that
\begin{equation*}
\sup \left\{ \left\Vert L\left( x\right) \right\Vert_{Y,X} \ :\ x\in B\right\}
<\frac{1-a'}{1-a}m_{0} \;. \label{aprime}
\end{equation*}
We define a continuous map $w:\, B' \to Y$ by the formula
$$w(y):=y-f\circ g_0(y) \in Y\;.$$
By continuity of $w$, the set
\[{\mathcal V}:= \left\{y\in \dot{B}'\ :\ \Vert w(y)\Vert_Y < \frac{1}{2} \varphi (g_0)\Vert y\Vert_Y\right\}\]
is open in $\dot{B}'\,.$
\medskip

Now, $Df$ is bounded since $f$ is Lipschitz continuous, and $L$ is bounded on $B$ by Assumption $(iii)$. Therefore, combining these bounds with the continuity of $w$, we see that for each $(x,y)\in B\times (\dot{B}'\setminus {\mathcal V})\,,$ there exists a positive radius $\beta(x,y)$ such that, if $(x',y')\in B_X(x,\beta(x,y))\times B_Y(y,\beta(x,y))\,,$
then $(x',y')\in B\times \dot{B}'$ and
\[\left(\Vert Df(x')\Vert_{X,Y}\Vert L(x)\Vert_{Y,X}+1+a'\right)\Vert w(y')-w(y)\Vert_X\leq (a'-a)\Vert w(y)\Vert_X\,,\]
which implies the inequality
\begin{equation}
a\Vert w(y)\Vert_Y+\left\Vert \left(Df(x')\circ L(x)-I_Y\right) (w(y')-w(y))\right\Vert_Y\leq a' \Vert w(y')\Vert_Y\,. \label{ytoyprime}
\end{equation}
Let $\gamma(x,y):=\min\left(\alpha(x,w(y)) ; \,\beta(x,y)\right)$ where $\alpha(x,w)$ is the radius introduced in Assumption $(ii)$. Then this assumption combined with (\ref{ytoyprime}) implies that 

\begin{equation}
\left\Vert \left( Df\left( x'\right) \circ L\left( x\right)-I_Y\right) w(y') \right\Vert_Y \leq a' \Vert w(y') \Vert_Y\label{8pourwprime}
\end{equation}
for each $(x,y)\in B\times (\dot{B}'\setminus {\mathcal V})$ and all $(x',y')\in B_X(x,\gamma(x,y))\times B_Y(y,\beta(x,y))\,.$
\medskip

Since the set $\Omega:=\bigcup_{(x,y)\in B\times (\dot{B}'\setminus {\mathcal V})}B_X(x,\gamma(x,y))\times B_Y(y,\beta(x,y))$ is a metric space, it is paracompact \cite{Stone}. Thus, $\Omega$ has a locally finite open covering $(\omega_i)_{i\in I}$ where for each $i\in I\,,$
$$\omega_i\subset B_X(x_i,\gamma(x_i,y_i))\times B_Y(y_i,\beta(x_i,y_i))$$
for some $(x_i,y_i)\in B\times (\dot{B}'\setminus {\mathcal V})\,.$ In the sequel, we take the norm $\max(\Vert x\Vert_X;\,\Vert y\Vert_Y)$ on $X\times Y$. For $(x,y)\in B\times \dot{B}'$ we define
\begin{eqnarray*}
\sigma_i(x,y)&:=&{\rm dist}((x,y),(B\times \dot{B}')\setminus \omega_i)\,, \\
\theta(x,y)&:=&\frac{\sum_{i\in I}\sigma_i(x,y)}{{\rm dist}(y,\dot{B}'\setminus {\mathcal V})+\sum_{i\in I}\sigma_i(x,y)}\in [0,1]\,.
\end{eqnarray*}
Note that $\theta(x,y)=1$ when $\Vert w(y)\Vert_Y\geq\frac{\varphi(g_0)}{2}\Vert y\Vert_Y\,$, and $\theta(x,y)=0$ when $(x,y)\notin \Omega\,.$\medskip

We are now ready to define
$$\tilde{L}(x,y)= \left({\rm dist}(y,\dot{B}'\setminus {\mathcal V})+\sum_{i\in I}\sigma_i(x,y)\right)^{-1}\sum_{i\in I}\sigma_i(x,y)L(x_{i})\,.$$
One easily checks that
$\tilde{L}$ is locally Lipschitz on $B\times \dot{B}'$. Moreover, it satisfies the same uniform estimate as $L$:
\begin{equation}
\sup \left\{ \Vert \tilde{L}\left( x,y\right) \Vert_{Y,X} \ :\ x\in B\,,\;y\in \dot{B}'\right\} <\frac{1-a'}{1-a}m_0
\label{uniform}
\end{equation}
and due to (\ref{8pourwprime}), it is an approximate right-inverse of $Df$ ``in the direction $w(y)$". More precisely, for all
$(x,y)$ in $B\times \dot{B}'$, one has:
\begin{equation}
\label{8bis}
\left\Vert \left( Df\left( x\right) \circ \tilde{L}\left( x,y\right)-\theta(x,y)I_Y\right) w(y)\right\Vert_Y
\leq a'\theta(x,y)\Vert w(y)\Vert_Y \;.
\end{equation}
Now, to each $y\in \dot{B}'$ we associate the vector field on $B$
$$X_y(x):=\tilde{L} ( x,y ) w(y)$$
and we consider the Cauchy problem
\begin{equation*}
\left \lbrace
\begin{aligned}
\frac{dx}{dt} &=&X_y(x),\\
x(0)&=&g_0(y).
\end{aligned}
\right.
\end{equation*}
The vector field $X_y$ is locally Lipschitz in the variable $x\in B$ and from \eqref{uniform} we have the uniform estimate
$$\sup \left\{ \Vert y\Vert_Y^{-1}\Vert X_y(x)\Vert_X \, :\, x\in B\,,\; y\in \dot{B}'\right\} < \frac{1-a'}{1-a}m_0  \, \varphi( g_{0}) \leq \frac{1-a'}{1-a}m_0 \;.$$
Thus, recalling that $\Vert g_0(y)\Vert_X < R m_0 m^{-1}\, ,$ we see that our Cauchy problem has a unique solution $x(t)=g_t(y)\in B$ on the time interval $[ 0,\tau]$ with $\tau=\frac{m-m_0 }{(1-a') m_0}\,.$ In addition, we take $g_t(0)=0$. This gives us a one-parameter family of functions $g_t:B'\to B$. For $0<t\leq \tau\,$, $g_t$ satisfies the estimate
\begin{equation} 
\sup \left\{ \Vert y\Vert_Y^{-1}\Vert g_t(y)-g_0(y)\Vert_X \, :\, y\in \dot{B}'\right\}  <  \frac{1-a'}{1-a}m_0  \, \varphi( g_{0})  t \;.\label{distance}
\end{equation}
Since $\left\Vert g_{0}\right\Vert_{\mathcal{C}} \leq \frac{m_0 }{1-a}$ and $\varphi(g_0)\leq 1$, the inequality \eqref{distance} implies that \begin{equation}\label{dom}
\sup\limits_{\dot{B}'} \Vert y\Vert_Y^{-1}\Vert g_t(y)\Vert_X <\frac{m}{1-a}\;,\ \forall t\in [0,\tau]\,.
\end{equation}
We recall that $\tilde{L}(\cdot,\cdot)$ is locally Lipschitz on $B\times \dot{B}'$ and $g_0$, $w$ are continuous.
Thus, by Gronwall's inequality, for each $t\in [0,\tau]\,$ the function $g_t$ is continuous on $\dot{B}'$. 
We can conclude that $g_t\in {\mathcal C}\,$, and \eqref{dom} implies that $\varphi(g_t)<\infty\,$.\medskip

Now, to each $(t,y)\in [0,\tau]\times \dot{B}'$ we associate $s_t(y):=\int_0^t \theta(g_u(y), y)du\,$ and we consider the function
\[
(t,y)\in [0,\tau]\times \dot{B}'\mapsto h(t,y):=f\circ g_t(y) -y +(1-s_t(y)) w(y)\in Y\,.
\]

Since $f$ is Lipschitzian, its G\^ateaux differential $Df(x)$ at any $x\in B$ is also a Hadamard differential, as mentioned in the introduction. This implies that for any function $\gamma:\,(-1,1)\to B$ differentiable at $0$ and such that $\gamma(0)=x$, the function $f\circ \gamma$ is differentiable at $0$
and the chain rule holds true: $(f\circ \gamma)'(0)=Df(\gamma(0))\gamma'(0)\,.$\medskip

Therefore, using (\ref{8bis}), we get:

\begin{eqnarray}
\left\Vert\frac{\partial}{\partial{t}}h(t,y)\right\Vert_Y 
&=&\Vert (Df\left( g_t(y) \right) \circ \tilde{L}\left( g_t(y),y  \right)-\theta(g_t(y),y)I_Y) w(y) \Vert_Y  \label{9} \\
&\leq &a'\theta(g_t(y),y)\Vert  w(y) \Vert_Y\;. \nonumber
\end{eqnarray}
In addition $h(0,y)=0$, so by the mean value theorem,
$$\Vert h(t,y)\Vert_Y \leq a's_t(y) \Vert  w(y) \Vert_Y\,.$$
By the triangle inequality, this implies that
\begin{equation}
\Vert f\left( g_t(y)\right) -y \Vert_Y \leq (1-(1-a')s_t(y)) \Vert  w(y) \Vert_Y\,.
\label{decay}
\end{equation}

We are now ready to get a contradiction. The estimate (\ref{distance}) may be written as follows: 
\begin{equation}
\frac{1-a}{m_0}\,\Vert g_t-g_0\Vert_{\mathcal C} <  (1-a') t \, \varphi( g_{0}) \;,\quad\forall t\in (0,\tau]\,.  \label{distancebis}
\end{equation}
As a first consequence of \eqref{decay}, if $\Vert y\Vert_Y^{-1}\Vert f\left( g_t(y)\right) -y \Vert_Y \geq \frac{\varphi(g_0)}{2}\,$ then one has $\,\Vert w(y)\Vert_Y\geq\frac{\varphi(g_0)}{2}\Vert y\Vert_Y\,$, hence $s_t(y)=t$. Thus, the estimate (\ref{decay}) implies that for all $(t,y)\in [0,\tau]\times \dot{B}'$,
\[
\Vert y\Vert_Y^{-1}\Vert f\left( g_t(y)\right) -y \Vert_Y \leq \max\left\{\frac{\varphi(g_0)}{2}\,;\,(1-(1-a')t)\Vert y\Vert_Y^{-1}\Vert w(y)\Vert_Y\right\} \,.
\]
Moreover we always have $\Vert y\Vert_Y^{-1}\Vert w(y)\Vert_Y\leq \varphi(g_0)$, so, with $\tau':=\min\left\{\tau,\frac{1}{2}\right\}$ we get
\[\Vert y\Vert_Y^{-1}\Vert f\left( g_t(y)\right) -y \Vert_Y \leq (1-(1-a')t) \varphi(g_0)\]
for all $0\leq t\leq \tau'$ and $y\in\dot{B}'$. This means that
\begin{equation}
\varphi(g_t) \leq (1-(1-a')t) \varphi(g_0)\;,\quad\forall t\in [0,\tau']\,.
\label{decaybis}
\end{equation}
Combining \eqref{distancebis} with \eqref{decaybis}, we find the following, for $0< t \leq \tau'\,$: 
$$\varphi(g_t)\leq \varphi (g_0)-(1-a')t\,\varphi(g_0) < \varphi (g_0)-\frac{1-a}{m_0}\,\Vert g_t-g_0\Vert_{\mathcal C}\,,$$
which contradicts (\ref{opti}). This ends the proof of the existence statement in Theorem \ref{continuous inverse}.\medskip

The uniqueness statement is proved by more standard arguments: if $g_1$ and $g_2$ are two continuous right-inverses of $f$ such that $g_1(0)=g_2(0)=0$, then the set
\[Z:=\{y\in B'\;|\;g_1(y)=g_2(y)\}\]
is nonempty and closed. On the other hand, if $Df(x)$ is left and right invertible, it is an isomorphism. By Remark \ref{inverse} its inverse $\hat{L}(x)$ is bounded independently of $x$. 
We fix an arbitrary $y_0$ in $Z$ and we consider a small radius $\rho>0$ (to be chosen later) such that $B_Y(y_0,\rho)\subset B'$. By continuity of $g_1-g_2$ at $y_0$, there is $\eta(\rho)>0$ such that $\lim\limits_{\rho\to 0}\eta(\rho)=0$ and, for each $y$ in the ball $B_Y(y_0,\rho)\,$,
\[
 \Vert g_2(y)-g_1(y)\Vert_X\leq \eta(\rho)\,.
\]
Moreover, we have $f(g_2(y))-f(g_1(y))=y-y=0$. Thus, using $(iv)$, we find that
\[
\Vert Df(g_1(y))(g_2(y)-g_1(y))\Vert_Y\leq (\varepsilon\circ\eta)(\rho)\,\Vert g_2(y)-g_1(y)\Vert_Y\,.
\]
Then, multiplying $Df(g_1(y))(g_2(y)-g_1(y))$ on the left by $\hat{L}(g_1(y))$ and using the uniform bound on $\hat{L}$, we get a bound of the form
\[
\Vert g_2(y)-g_1(y)\Vert_X\leq \xi(\rho)\,\Vert g_2(y)-g_1(y)\Vert_X
\]
with $\lim\limits_{\rho\to 0}\xi(\rho)=0$. As a consequence, for $\rho$ small enough one has $g_2(y)-g_1(y)=0$, so $y\in Z$.
This proves that $Z$ is open. By connectedness of $B'$ we conclude that $Z=B'$, so $g_1$ and $g_2$ are equal. This ends the proof of Theorem \ref{continuous inverse}.

\section{A hard surjection theorem with continuous right-inverse}\label{hard} In this section we state our hard surjection theorem with continuous right-inverse and we shortly explain its proof which is a variant of the arguments of \cite{ES} in which Theorem \ref{Ivar} is replaced by Theorem \ref{continuous inverse}.

Let $(V_{s},\,\Vert\cdot\Vert_{s})_{0\leq s \leq S}$ be a scale of
Banach spaces, namely:
\[
0\leq s_{1}\leq s_{2}\leq S\Longrightarrow\left[  V_{s_{2}}\subset
V_{s_{1}}\text{ \ and\ \ }\Vert\cdot\Vert_{s_{1}}\leq\Vert\cdot\Vert_{s_{2}%
}\right]\;.
\]

We shall assume that to each $\Lambda\in [1,\infty)$ is associated a continuous linear projection $\Pi(\Lambda)$ on
$V_0$, with range $E(\Lambda)\subset V_S$. We shall also assume that the spaces $E(\Lambda)$ form a nondecreasing family of sets indexed by $[1,\infty)$, while the spaces $Ker\,\Pi(\Lambda)$ form a nonincreasing family. In other words:

\[1\leq \Lambda\leq \Lambda'\,\Longrightarrow \,\Pi(\Lambda)\Pi(\Lambda')=\Pi(\Lambda')\Pi(\Lambda)=\Pi(\Lambda)\;.\]

Finally, we assume that the projections $\Pi(\Lambda)$ are ``smoothing operators" satisfying the following estimates:
\medskip

\textbf{Polynomial growth and approximation}: \textit{There are constants }$A_{1},\ A_{2}\geq
1$\textit{ such that, for all numbers} $0\leq s\leq S$\textit{, all }$\Lambda\in [1,\infty)$\textit{ and all }$u\in V_{s}\,$\textit{, we have:}%
\begin{align}
\forall t\in [0,S]\,,\;\;\Vert\Pi(\Lambda)u\Vert_{t}  &  \leq A_{1}\,\Lambda^{(t-s)^{+}}\Vert u\Vert_{s}\,.
\label{loss}\\
\forall t\in [0,s]\,,\;\;\Vert(1-\Pi({\Lambda}))u\Vert_{t}  &  \leq A_{2}\,\Lambda^{-(s-t)}\Vert u\Vert_{s}\,.
\label{gain}
\end{align}

When the above properties are met, we shall say that $(V_{s}\,,\,\Vert\cdot\Vert_{s})_{0\leq s \leq S}$
endowed with the family of projectors
$\left\{\,\Pi(\Lambda)\;,\;\Lambda\in [1,\infty)\,\right\}\,, $ is a {\it tame} Banach scale.
\medskip

Let $(W_{s}\,,\,\Vert\cdot\Vert'_s)_{0\leq s\leq S}$ be another tame scale of Banach spaces. We shall denote by $\Pi^{\prime}(\Lambda)$ the corresponding projections defined on $W_0$ with ranges $E'(\Lambda)\subset W_S$, and by $A'_{i}\;(i=1,2,3)$ the
corresponding constants in (\ref{loss}), (\ref{gain}).\medskip

We also denote by $B_{s}$ the unit ball in
$V_{s}$ and by $B'_s(0,r)$ the ball of center $0$ and positive radius $r$ in $W_s$:
\[
B_{s}=\left\{  u\in V_s\ | \ \left\Vert u\right\Vert _{s}< 1\right\}\quad \hbox{ and } \quad B'_s(0,r)=\left\{  v\in W_s \ | \ \left\Vert v\right\Vert' _{s}< r\right\}\;.
\]

In the sequel we fix nonnegative constants $s_0, m,\ell$ and $\ell'$. We will assume that $S$ is large enough.\bigskip

We first recall the definition of G\^ateaux-differentiability, in a form adapted to our framework:\medskip

\begin{definition}\label{Gdiff}
We shall say that a function $F:\,B_{s_{0}+m}\rightarrow
W_{s_{0}}$ is \emph{G\^{a}teaux-differentiable} (henceforth
G-differentiable) if for every $u\in
B_{s_{0}+m}$, there exists a linear map $DF\left(
u\right)  :V_{s_0+m}\rightarrow W_{s_0}$ such that for every
$s\in [s_0,S-m]$, if $u\in
B_{s_{0}+m}\cap V_{s+m}$,  then $DF \left(
u\right)$ maps continuously $V_{s+m}$ into $W_s$, and %
\[
\forall h\in V_{s+m}\ ,\ \lim_{t\rightarrow0}\ \left\Vert \frac{1}{t}\left[
F\left(  u+th\right)  -F\left(  u\right)  \right]  -D F\left(  u\right)
h\right\Vert _{s}^{\prime}=0\;.
\]

\end{definition}

\medskip
Note that, even in finite dimension, a G-differentiable map need not be $C^{1}$, or
even continuous. However, if $D F:\,B_{s_{0}+m}\cap V_{s+m}\to\mathcal{L}(V_{s+m},W_s)$ is locally bounded, then $F:\,B_{s_{0}+m}\cap V_{s+m}\to W_s$ is locally Lipschitz, hence continuous. In the present paper, we are in such a situation.\bigskip

We now define the notion of $S$-tame differentiability:\bigskip

\begin{definition}
\leavevmode\par
\begin{itemize}

\item We shall say that the map $F:\,B_{s_{0}+m}\rightarrow
W_{s_{0}}\,$ is $S$-tame differentiable if it is G-differentiable in the sense of Definition \ref{Gdiff}, and, for some positive constant $a\,$ and all $s \in [s_{0}, S-m]\,$, if $u\in B_{s_{0}+m}\cap V_{s+m}$ and $h\in V_{s+m}\,,$ then
$DF\left(  u\right)  h \in W_s$ with the tame direct estimate
\begin{equation}
\left\Vert DF\left(  u\right)  h\right\Vert _{s}^{\prime}\leq
a\left(  \left\Vert h\right\Vert _{s+m}+\left\Vert u\right\Vert _{s+m}
\left\Vert h\right\Vert _{s_{0}+m}\right)\;.  \label{tamedirect}
\end{equation}

\item Then we shall say that $DF$ is tame right-invertible if there are $\,b>0$ and $\ell,\,\ell'\geq 0$ such that for all
$u\in B_{s_{0}+\max\{m,\ell\}}$ , there is a linear map
$L\left(  u\right)  :W_{s_0+\ell^{\prime}}\rightarrow V_{s_0}$
satisfying
\begin{equation}
\forall k\in W_{s_0+\ell'}\,,\ \ \ DF
\left(  u\right)
L\left(  u\right)  k=k \label{rightinverse}
\end{equation}
and for all $s_{0}\leq s\leq
S-\max\left\{  \ell,\ell^{\prime}\right\}  $, if $u\in B_{s_{0}+\max\{m,\ell\}}\cap V_{s+\ell}\,$ and $k\in W_{s+\ell^{\prime}}\,,$ then $L\left(
u\right)  k\in V_s\,,$ with the tame inverse estimate
\begin{equation}
\left\Vert L\left(
u\right)  k\right\Vert _{s}\leq b\left(  \left\Vert
k\right\Vert _{s+\ell^{\prime}}^{\prime}+\left\Vert k\right\Vert _{s_{0}
+\ell^{\prime}}^{\prime}\left\Vert u\right\Vert _{s+\ell}\right)\;.
\label{tameinverse}
\end{equation}

\end{itemize}
\end{definition}
\medskip

In the above definition, the numbers $m, \ell, \ell^{\prime}$ represent the loss of derivatives for $DF$ and its right-inverse.\medskip

The main result of this section is the following theorem.

\begin{theorem}
\label{Thm8}

Assume that the map $F:\,B_{s_{0}+m}\rightarrow
W_{s_{0}}\,$ is $S$-tame differentiable between the tame scales $(V_s)_{0\leq s\leq S}$ and $(W_s)_{0\leq s\leq S}$ with $F(0)=0$ and that $DF$ is tame right-invertible. Let $s_0,\,m,\,\ell,\,\ell'$ be the associated parameters.\medskip

Assume in addition that for each $\Lambda,\Lambda'\in [1,S]$ the map \[u\in B_{s_{0}+\max\{m,\ell\}}\cap E(\Lambda)\mapsto \Pi'_{\Lambda'} DF(u)\!\upharpoonright_{_{E_\Lambda}}\in\mathcal L(E(\Lambda),E'(\Lambda'))\]
is continuous for the norms $\Vert\cdot\Vert_{s_0}$ and $\Vert\cdot\Vert'_{s_0}$.\medskip

Let $s_1 \geq s_0+\max\{m,\ell\}$ and $\delta>s_1+\ell'$.
Then, for $S$ large enough, there exist a radius $r>0$ and a continuous map
$G: B'_\delta(0,r)\to B_{s_1}$ such that:
\begin{align*}
&\;G(0)=0\quad \hbox{and}\quad F\circ G=I_{B'_\delta(0,r)}\,. \\
&\left\Vert G(v)\right\Vert _{s_1}\leq r^{-1}\,\left\Vert v\right\Vert _{\delta}^{\prime}\;,\quad \forall v\in B'_\delta(0,r)\,.
\end{align*}

\end{theorem}

As mentioned in the introduction, compared with the results of \cite{ES} the novelty in Theorem \ref{Thm8} is the continuity of $G$. To prove this theorem, one repeats with some modifications the arguments of \cite{ES} in the case $\varepsilon=1$ (in that paper a singularly perturbed problem depending on a parameter $\varepsilon$ was dealt with, but for simplicity we do not consider such a dependence here). With the notation of that paper, let us explain briefly the necessary changes.
\medskip

We recall that in \cite{ES} a vector $v$ was given in $B'_\delta(0,r)$ and the goal was to solve the equation $F(u)=v$. The solution $u$ was the limit of a sequence $u_n$ of approximate solutions constructed inductively. Each $u_n$ was a solution of the projected equation $\Pi'_nF(u_{n})=\Pi'_{n-1}v$, $\,u_n\in E_{n}$. It was found as $u_n=u_{n-1}+z_n$, $\,z_n$ being a small solution in $E_n$ of an equation of the form $f_n(z)=\Delta_n v + e_n$, with $f_{n}\left(z\right):=\Pi'_{n}\left(F\left(u_{n-1}+z\right)-F\left(u_{n-1}\right)\right)\,$,
$\ \Delta_n v:=\Pi'_{n-1}(1-\Pi'_{n-2})v\,$ and $e_n:=-\Pi'_{n}(1-\Pi'_{n-1})F(u_{n-1})\,$.
The existence of $z_n$ was proved by applying Theorem \ref{Ivar} to the function $f_n$ in a ball $B_{\mathcal N_n}(0, R_n)$ (see \cite{ES}[\S 3.3] for precise definitions of $E_n$, $\Pi'_n$, $f_n$ and $\mathcal N_n$).\medskip

Instead, in order to prove Theorem \ref{Thm8} we construct inductively a sequence of continuous functions $G_n: B'_{\delta}(0,r)\to B_{s_1}\cap E_n $ such that $\Pi'_n F\circ G_n(v)=\Pi'_{n-1}v$ for all $v$ in $B'_{\delta}(0,r)$. Each $G_n$ is of the form $G_{n-1}+H_n$ with
\[H_n(v)=g_n\left(\Delta_n v -\Pi'_n(1-\Pi'_{n-1})F\circ G_{n-1}(v)\right)\]
where $g_n$ is a continuous right-inverse of $f_n$ such that $g_n(0)=0$, obtained thanks to Theorem \ref{continuous inverse}.\medskip

Under the same conditions on the parameters as in \cite{ES}, we find that the sequence of continuous functions $(G_n)_n$ converges uniformly on $B'_{\delta}(0,r)$ for the norm $\Vert\cdot\Vert_{s_1}$ and this implies the continuity of their limit $G:\,B'_{\delta}(0,r)\to B_{s_1}$. This limit is the desired continuous right inverse of $F$. We insist on the fact that the conditions on $r$ are exactly the same as in \cite{ES}. Indeed, in order to apply Theorem \ref{continuous inverse} to $f_n$ we just have to check assumptions $(i)$, $(ii')$ and $(iii)$. This is done with exactly the same constraints on the parameters as in \cite{ES}. $\square$
\bigskip\bigskip

We end the paper with a uniqueness result, which requires additional conditions.

\begin{theorem}
\label{Thm9}

Suppose that we are under the assumptions of Theorem \ref{Thm8}, and that the following two additional conditions hold true:\medskip

$\bullet$ For each $u\in B_{s_0+\max(m,\ell)}\,$:
\begin{equation}
\forall h\in V_{s_0+m+\ell'}\,,\ \ \ L\left(  u\right) DF
\left(  u\right)
  h=h\,. \label{leftinverse}
\end{equation}

$\bullet$ For each $s\in [s_0,S-m]$ and each $c>0$ there exists a non-decreasing function $\varepsilon_{s,c}:\,(0,\infty)\to (0,\infty)$ such that $\lim\limits_{t\to 0}\varepsilon_{s,c}(t)=0$ and, for all $u_1,\,u_2$ in  $B_{s_0+m}\cap E_{s+m}\,$ with $\Vert u_1\Vert_{s+m}\leq c\,$:
\begin{equation}\label{tameder}
\Vert F(u_2)-F(u_1)-DF(u_1)(u_2-u_1)\Vert_s\leq \varepsilon_{s,c}(\Vert u_2-u_1\Vert_{s+m}) \Vert u_2-u_1\Vert_{s_0+m}\,.
\end{equation}
Let $s_1 \geq s_0+\max\{2m+\ell',m+\ell\}$.
Then, for any $S\geq s_1$, $\delta\in [s_0,S]$ and $r>0$, there is at most one map
$G: B'_\delta(0,r)\to B_{s_0+\max(m,\ell)}\cap W_{s_1}$ continuous for the norms $\Vert\cdot\Vert'_{\delta}$ and $\Vert\cdot\Vert_{s_1}\,$, such that
\begin{equation}
\label{cond-right}
G(0)=0\quad \hbox{and}\quad F\circ G=I_{B'_\delta(0,r)}\,. 
\end{equation}

\end{theorem}

\begin{remark} The tame estimate \eqref{tameder} is satisfied, in particular, when $F$ is of class $C^2$ with a classical tame estimate on its second derivative as in \cite[(2.11)]{TZ}. In that special case, for $s$ and $c$ fixed one has the bound $\varepsilon_{s,c}(t)=O(t)_{t\to 0}\,$.
\end{remark}

In order to prove Theorem \ref{Thm9}, we assume that $G_1$, $G_2$ both satisfy \eqref{cond-right} and we introduce the set $Z:=\{v\in B'_\delta(0,r)\;|\; G_1(v)=G_2(v)\}$. This set is nonempty since it contains $0$, and it is closed in $B'_\delta(0,r)$ for the norm $\Vert\cdot\Vert'_{\delta}$ by continuity of $G_1-G_2$. It remains to prove that it is open.\medskip

For that purpose, we fix an arbitrary $v_0$ in $Z$ and we consider a small radius $\rho>0$ (to be chosen later) such that $B'_\delta(v_0,\rho)\subset B'_\delta(0,r)$. By continuity of $G_1\,$, $G_2$ at $v_0$, there is $\eta(\rho)>0$ such that $\lim\limits_{\rho\to 0}\eta(\rho)=0$ and, for each $v$ in the ball $B'_\delta(v_0,\rho)\,$:
\[
 \Vert G_1(v)\Vert_{s_1}\leq \Vert G_1(v_0)\Vert_{s_1}+\eta(\rho)\ \hbox{ and }\ \Vert G_2(v)-G_1(v)\Vert_{s_1}\leq \eta(\rho)\,.
\]
However, we also have $F(G_2(v))-F(G_1(v))=v-v=0$. Thus, imposing $\eta(\rho)\leq 1$ and applying \eqref{tameder} with $s=s_1-m$, $c=\Vert G_1(v_0)\Vert_{s_1}+1$ and $u_i=G_i(v)\,$, $i=1,2$, we find that
\[
\Vert DF(G_1(v))(G_2(v)-G_1(v))\Vert'_{s_1-m}\leq (\varepsilon_{s_1-m,c}\circ\eta)(\rho)\,\Vert G_2(v)-G_1(v)\Vert_{s_0+m}\,.
\]
Then, multiplying $DF(G_1(v))(G_2(v)-G_1(v))$ on the left by $L(G_1(v))$ and using \eqref{leftinverse} and the tame estimate \eqref{tameinverse}, we get a bound of the form
\[
\Vert G_2(v)-G_1(v)\Vert_{s_1-\max(m+\ell',\ell)}\leq \xi(\rho)\,\Vert G_2(v)-G_1(v)\Vert_{s_0+m}
\]
with $\lim\limits_{\rho\to 0}\xi(\rho)=0$. Since $s_1-\max(m+\ell',\ell)\geq s_0+m$, we conclude that for $\rho$ small enough one has $G_2(v)-G_1(v)=0$, so $v\in Z$. The set $Z$ is thus nonempty, closed and open in $B'_\delta(0,r)$, so we conclude that $Z=B'_\delta(0,r)$ and Theorem \ref{Thm9} is proved. $\square$
\bigskip\bigskip

\end{document}